             \documentclass[11pt]{article}
\usepackage{hyperref,amsmath,amsthm,amsfonts,amssymb,amscd,enumerate,tikz}

\usetikzlibrary{positioning}

\theoremstyle{plain}

\newtheorem{corollary}{Corollary}

\newtheorem*{Theorem}{Theorem}

\DeclareMathOperator{\piprod}{\raisebox{-0.1em}{\huge{$\pi$}}\kern -0.2em}

\theoremstyle{definition}

\newtheorem*{example*}{Example}

\newtheorem*{Remarks}{Remarks}

\newcommand{\cg}{\mathcal {G}}

\newcommand{\cp}{\mathcal {P}}

\newcommand{\cs}{\mathcal {S}}

\newcommand{\cz}{\mathcal {Z}}

\newcommand{\bA}{{\mathbf A}}
\newcommand{\bB}{{\mathbf B}}

\newcommand{\bx}{{\mathbf x}}

\newcommand{\bG}{{\mathbf G}}

\newcommand{\be}{{\mathbf e}}

\newcommand{\bq}{{\mathbf q}}

\newcommand{\bs}{{\mathbf s}}
\newcommand{\bt}{{\mathbf t}}

\newcommand{\m}{{[m]}}
\newcommand{\ah}{\hat{h}}

\newcommand{\bone}{{\mathbf 1}}

\newcommand{\zl}{{\cz_L}}

\newcommand{\zx}{{\zl(\bA,\bB)}}

\newcommand{\ab}{(\bA,\bB)}

\def\clap#1{\hbox to 0pt{\hss#1\hss}}

\newcommand{\comment}[1]{}

\newcommand{\cone}{\operatorname{Cone}}

\newcommand{\supp}{\operatorname{Supp}}

\newcommand{\Econe}{(\boldsymbol{\cone}\, \mathbf E, \mathbf E)}

\newcommand{\bman}{(\mathbf M, \boldsymbol{\partial}\mathbf M)}
\newcommand{\dman}{\boldsymbol{\partial}\mathbf M}

	{
\end{enumerate}
}
\newenvironment{enumerate1*}{
\begin{enumerate}[\upshape (*1)]}%
	{
\end{enumerate}
}
	{
\end{enumerate}
}
\newenvironment{enumeratei'}{
\begin{enumerate}[\upshape (i)$'$]}%
	{
\end{enumerate}
}

\newenvironment{enumeratea'}{
\begin{enumerate}[\upshape (a)$'$]}{
\end{enumerate}
}

\begin{document}
\title{The Euler characteristic of a polyhedral product}
\author{Michael W. Davis\thanks{The  author was partially supported by NSF grant DMS-1007068 and the Institute for Advanced Study.}  }
\date{\today} \maketitle
\begin{abstract}
Given a finite simplicial complex $L$ and a collection of pairs of spaces indexed by its vertex set, one can define their polyhedral product.  We record a  simple formula for its Euler characteristic.  In special cases the formula simplifies further to one involving the $h$-polynomial of $L$.
	 \smallskip

	\paragraph{AMS classification numbers.} Primary:  20F65, 57M07. \\
	Secondary:  13F55.

	\paragraph{Keywords:}  Charney-Davis Conjecture, Euler characteristic, graph product, polyhedral product. 
\end{abstract}

The purpose of this note is to record a simple formula for the Euler characteristic of the polyhedral product of a collection of pairs of spaces.  
To define this notion, start with a finite simplicial complex $L$ with vertex set $\m$, where $[m]:=\{1,\dots, m\}$.  Let $\cs(L)$ be the poset (of vertex sets) of simplices in $L$, including the empty simplex.  Next suppose we are given a family of pairs of finite CW complexes, $(\bA,\bB)=\{(A_i,B_i)\}_{i\in \m}$.   Denote the product $\prod_{i\in [m]}A_i$ by $\prod \bA$.  For $\bx:=(x_i)_{i\in I}$ a point in $\prod \bA$, put 
	\(
	\supp(\bx):=\{i\in \m\mid x_i\notin B_i\}.  
	\)
Define the \emph{polyhedral product}, $\zx$, by 
	\[
	\zx=\{\bx\in\prod\bA\mid \supp(\bx)\in \cs(L)\}.
	\]
(The terminology comes from \cite{bbcg2}.) 
Define $\m$-tuples $\be(\bA)$ and $\be(\bB)$ by
	\[
	\be(\bA)=(\chi(A_i))_{i\in \m} \quad\text{and} \quad \be(\bB)=(\chi(B_i))_{i\in \m}
	\]
	
Let  $\cp(I)$ denote the power set of a finite set $I$. Given an $I$-tuple $\bt=(t_i)_{i\in I}$ and $J\in \cp(I)$, define a monomial $\bt_J$ by
	\[
	\bt_J=\prod_{j\in J} t_j
	\]
The following is a version of the Binomial Theorem,
	\begin{equation}\label{e:binom}
	(\bs+\bt)_I= \sum_{J\in \cp(I)} \bs_J \bt_{I-J}.
	\end{equation}	

\begin{Theorem}
\[
\chi (\zl\ab)= \sum_{J\in \cs(L)}(\be(\bA)-\be(\bB))_J\be(\bB)_{\m-J},
\]
\end{Theorem}

\begin{proof}
For $I\in \cs(L)$, the Binomial Theorem \eqref{e:binom} gives
\[
\chi(\prod_{i\in I} A_i)=\be(\bA)_I=((\be(\bA)-\be(\bB)) +\be(\bB))_I
=\sum_{J\in\cp(I)} (\be(\bA)-\be(\bB))_J\be(\bB)_{I-J}.
\]
In other words, the contribution to the Euler characteristic  from the part of $\zl\ab$ corresponding to the open simplex on $I$ is $(\be(\bA)-\be(\bB))_I\be(\bB)_{\m-I}$. The formula follows.
\end{proof}

The \emph{$f$-polynomial} of $L$ is the polynomial in indeterminates $\bt=(t_i)_{i\in \m}$ defined by
	\begin{equation*}
	f_L(\bt)=\sum_{J\in \cs(L)} \bt_J .
	\end{equation*}
The \emph{$\ah$-polynomial} of $L$ is defined by
	\[
	\ah_L(\bt):=(\bone-\bt)_\m \,f_L\left(\frac{\bt}{\bone -\bt}\right).
	\]
If $\bt$ is the constant indeterminate given by $t_i=t$, then $f_L$ is a polynomial in one variable.   Denote it $f_L(t)$.  For $d=\dim L+1$, the \emph{$h$-polynomial} is defined by  $h_L(t)= \ah_L(t)/(1-t)^{m-d}=(1-t)^d f_L(\frac{t}{1-t})$.

\begin{corollary}\label{c:point}
Suppose each $B_i$ is a point $*_i$.  Then
\[
\chi (\zl(\bA,*))= \sum_{J\in \cs(L)}\tilde{\be}(\bA)_J =f_L(\tilde{\be}(\bA)),
\]
where $\tilde{\be}(\bA):=\be(\bA)-\bone=(\chi(A_i)-1)_{i\in \m}$ is the $m$-tuple of reduced Euler characteristics. 
\end{corollary}

\begin{corollary}\label{c:cone}
For each $i\in \m$, suppose $B_i=E_i$ a finite set of cardinality $q_i+1$ and $A_i=\cone E_i$. Put $\Econe=\{(\cone E_i, E_i)\}_{i\in \m}$. Then $\chi(\zl\Econe)=\ah_L(-\bq)$.
\end{corollary}

\begin{proof}
Any cone is contractible and hence, has Euler characteristic 1.  So, by the theorem,
	\begin{align*}
	\chi(\zl\Econe)&=\sum_{J\in \cs(L)} (-1)^{|J|} \bq_J (\bone+\bq)_{\m-J}\\
	&=(\bone+\bq)_\m\, f_L\left(\frac{-\bq}{\bone+\bq}\right) =\ah_L(-\bq).
	\end{align*}
\end{proof}

\begin{corollary}\label{c:mfld}
For each $i\in \m$, suppose $A_i$ is an odd-dimensional manifold with nonempty boundary $B_i$.  Then, for $d=\dim L +1$,
\[
\chi(\zl\ab)=\be(\bB)_\m \frac{h_L(-1)}{(-2)^d}.
\]
\end{corollary}

\begin{proof}
By Poincar\'e duality, $\chi(B_i)=2\chi(A_i)$. So, $\be(\bA)-\be(\bB)=-\frac{1}{2} \be(\bB)$.  By the theorem,
	\begin{align*}
	\chi (\zl\ab)&= \sum_{J\in \cs(L)} (-1/2)^{|J|}\be(\bB)_\m\\
	&=\be(\bB)_\m f_L(-1/2) =\be(\bB)_\m \,\frac{h_L(-1)}{(-2)^d}.
	\end{align*}
	\end{proof}

\begin{Remarks}
Here are some applications of these corollaries to aspherical spaces.  The basic fact is proved in \cite{d11}: if the following three conditions hold, then $\zl\ab$ is aspherical.
	\begin{itemize}
	\item
	For each $i \in \m$, $A_i$ is aspherical.
	\item
	For each $i\in \m$, each path component of $B_i$ is aspherical and any such component maps $\pi_1$-injectively into $A_i$.
	\item
	$L$ is a flag complex.
	\end{itemize}
When these conditions hold, $\zx$ is the classifying space of a group $\cg$, called the  ``generalized graph product'' of the $\pi_1(A_i)$ (cf.\ \cite{d11}).  In the next two paragraphs $\bG=(G_i)_{i\in \m}$ is an $m$-tuple of discrete groups and $L$ is a flag complex.

1) Suppose $\cg$ is the graph product of the $G_i$ with respect to the $1$-skeleton of $L$.  Then $B\cg$ is the polyhedral product of the classifying spaces $BG_i$ of the $G_i$ (cf.\ \cite{densuc}, \cite{do10}, \cite{d11}).  Suppose each $G_i$ is type FL (so that its Euler characteristic is defined).  Applying Corollary~\ref{c:point} to the case $(\bA,*)=\{(BG_i,*_i)\}_{i\in \m}$, we get
\[
\chi(B\cg)=f_L(\tilde{\be}(\bG)).
\]

2) Suppose each $G_i$ is finite of order $q_i+1$.  Let $\cg_0$ be the kernel of the natural epimorphism $\cg\to \prod_{i\in\m} G_i$.  Then $B\cg_0=\zl\Econe$, where $E_i=G_i$.  Applying Corollary~\ref{c:cone}, we get the following formula for the rational Euler characteristic of $\cg$ (cf.\ \cite[p.~308]{dbook}),
\[
\chi(\cg)=\frac{\chi(B\cg_0)}{(\bone+\bq)_\m}=\frac{\ah_L(-\bq)}{(\bone+\bq)_\m}.
\]

3)  The Euler Characteristic Conjecture asserts that if $N^{2k}$ is a closed aspherical $2k$-manifold, then $(-1)^k\chi(N^{2k})\ge 0$ (cf.\ \cite[p.~310]{dbook}).  The Charney-Davis Conjecture (cf.~\cite[p.~313]{dbook} asserts that if $L$ is a flag triangulation of a $(2c-1)$-sphere or even a ``generalized homology sphere'' (as defined in \cite[p.~192]{dbook}), then $(-1)^c h_L(-1)\ge 0$.  The point here is that if the Charney-Davis Conjecture is true, then one cannot find a counterexample to the Euler Characteristic Conjecture by using the construction in Corollary~\ref{c:mfld}.
Indeed, suppose $A_i=M_i$, a $(2k_i+1)$-manifold with (nonempty) boundary and $B_i=\partial M_i$.  Write $\bman$ for $\{(M_i,\partial M_i)\}_{i\in \m}$.   By Corollary~\ref{c:mfld},
	\begin{equation}\label{e:bman}
	\chi(\zl\bman)=\be(\dman)_\m \frac{h_L(-1)}{(-2)^d}.
	\end{equation}
If $L$ is a triangulation of a $(d-1)$-sphere (or generalized homology sphere and at least one $k_i$ is positive), then $\zl\bman$ is a closed manifold of dimension
	\[
	d+ \sum_{i=1}^m 2k_i .
	\]
This is an even integer if and only if $d$ is even; so, let us assume $d=2c$. Suppose, in addition, $\bman$ and $L$ satisfy the three conditions at the beginning of these remarks, so that $\zl\bman$ is a closed aspherical manifold.  Suppose that each $\partial M_i$ satisfies the Euler Characteristic Conjecture, i.e., $(-1)^{k_i}\chi(\partial M_i)\ge 0$.  So, the sign of $\be(\dman)_\m$ is $(-1)^{\sum k_i}$.  The Euler Characteristic Conjecture for $\zl\bman$ is 
\[
(-1)^{c+\sum k_i} \chi(\zl\bman) \ge 0.
\]
By \eqref{e:bman}, the sign of the left hand side is the sign of $(-1)^c h_L(-1)$.  So, the Euler Characteristic Conjecture for $\zl\bman$ is equivalent to $(-1)^ch_L(-1)\ge 0$, which is the Charney-Davis Conjecture.
\end{Remarks}

\end{document}